\documentclass[12pt,a4paper]{amsart}
\usepackage{amssymb}
\usepackage{amsmath}
\allowdisplaybreaks[2] 
\makeatletter
\let\Thebibliography=\thebibliography
\renewcommand{\thebibliography}[1]{\def\@mkboth##1##2{}\Thebibliography{#1}
\addcontentsline{toc}{section}{References}
\frenchspacing 
\setlength{\@topsep}{0pt}
\setlength{\itemsep}{0pt}%
\setlength{\parskip}{0pt plus 2pt}%
}
\makeatother
\makeatletter
\def\mdots@{\mathinner.\nonscript\!.%
 \ifx\next,.\else\ifx\next;.\else\ifx\next..\else
 \nonscript\!\mathinner.\fi\fi\fi}
\let\ldots\mdots@
\let\cdots\mdots@
\let\dotso\mdots@
\let\dotsb\mdots@
\let\dotsm\mdots@
\let\dotsc\mdots@
\def\vdots{\vbox{\baselineskip2.8\p@ \lineskiplimit\z@
    \kern6\p@\hbox{.}\hbox{.}\hbox{.}\kern3\p@}}
\def\ddots{\mathinner{\mkern1mu\raise8.6\p@\vbox{\kern7\p@\hbox{.}}%
    \raise5.8\p@\hbox{.}\raise3\p@\hbox{.}\mkern1mu}}
\makeatother
\makeatletter
\let\Enumerate=\enumerate
\renewcommand{\enumerate}{\Enumerate%
\setlength{\@topsep}{0pt}
\setlength{\itemsep}{0pt}%
\setlength{\parskip}{0pt plus 1pt}%
\renewcommand{\theenumi}{\textup{(\alph{enumi})}}%
\renewcommand{\labelenumi}{\theenumi}%
}
\let\endEnumerate=\endenumerate
\renewcommand{\endenumerate}{\endEnumerate\unskip}
\makeatother
\makeatletter
\def\@seccntformat#1{\csname the#1\endcsname.\quad}
\makeatother
\newcommand{\art}[6]{{\sc #1, \rm #2, \it #3 \bf #4 \rm (#5), \mbox{#6}.}}
\newcommand{\book}[3]{{\sc #1, \it #2, \rm #3.}}
\newcommand{\AND}{{\rm and }}
\RequirePackage{amsthm}
\newtheoremstyle{descriptive}%
  {\topsep}   
  {\topsep}   
  {\rmfamily} 
  {}          
  {\bfseries} 
  {.}         
  { }         
  {}          
\newtheoremstyle{propositional}%
  {\topsep}   
  {\topsep}   
  {\itshape}  
  {}          
  {\bfseries} 
  {.}         
  { }         
  {}          
\theoremstyle{propositional}
\newtheorem{thm}{Theorem}[section]
\newtheorem{prop}[thm]{Proposition}

\newtheorem{lemma}[thm]{Lemma} 
\newtheorem{cor}[thm]{Corollary}
\newtheorem{corollary}[thm]{Corollary} 
\theoremstyle{descriptive}
\newtheorem{deff}[thm]{Definition}

\newtheorem{remark}[thm]{Remark}
\makeatletter
\renewenvironment{proof}[1][\proofname]{\par
  \pushQED{\qed}%
  \normalfont
  \trivlist
  \item[\hskip\labelsep
        \itshape
    #1\@addpunct{.}]\ignorespaces
}{%
  \popQED\endtrivlist\@endpefalse
}
\makeatother
\newdimen\extrawidth
\def\iintlim#1#2{\setbox0\hbox{$\scriptstyle#1$}%
        \setbox1\hbox{$\scriptstyle#2$}%
        \extrawidth=\wd1 \advance\extrawidth-\wd0
        \ifdim\extrawidth<0pt \extrawidth=0pt\fi%
        \int_{#1\kern\extrawidth \kern .5em}^{#2\kern -\wd1} \kern -.5em%
}
\newcommand{\setm}{\setminus}
\def\vint_#1{\mathchoice%
          {\mathop{\kern 0.2em\vrule width 0.6em height 0.69678ex depth -0.58065ex
                  \kern -0.8em \intop}\nolimits_{\kern -0.4em#1}}%
          {\mathop{\kern 0.1em\vrule width 0.5em height 0.69678ex depth -0.60387ex
                  \kern -0.6em \intop}\nolimits_{#1}}%
          {\mathop{\kern 0.1em\vrule width 0.5em height 0.69678ex depth -0.60387ex
                  \kern -0.6em \intop}\nolimits_{#1}}%
          {\mathop{\kern 0.1em\vrule width 0.5em height 0.69678ex depth -0.60387ex
                  \kern -0.6em \intop}\nolimits_{#1}}}
\def\vintslides_#1{\mathchoice%
          {\mathop{\kern 0.1em\vrule width 0.5em height 0.697ex depth -0.581ex
                  \kern -0.6em \intop}\nolimits_{\kern -0.4em#1}}%
          {\mathop{\kern 0.1em\vrule width 0.3em height 0.697ex depth -0.604ex
                  \kern -0.4em \intop}\nolimits_{#1}}%
          {\mathop{\kern 0.1em\vrule width 0.3em height 0.697ex depth -0.604ex
                  \kern -0.4em \intop}\nolimits_{#1}}%
          {\mathop{\kern 0.1em\vrule width 0.3em height 0.697ex depth -0.604ex
                  \kern -0.4em \intop}\nolimits_{#1}}}
\DeclareMathOperator{\diam}{diam}

\DeclareMathOperator{\Lip}{Lip}

\DeclareMathOperator{\spt}{supp}
\DeclareMathOperator{\LLC}{LLC}

\DeclareMathOperator{\Capc}{Cap}

\newcommand{\loc}{_{\rm loc}}
{\catcode`p =12 \catcode`t =12 \gdef\eeaa#1pt{#1}}      
\def\accentadjtext#1{\setbox0\hbox{$#1$}\kern   
                \expandafter\eeaa\the\fontdimen1\textfont1 \ht0 }
\def\accentadjscript#1{\setbox0\hbox{$#1$}\kern 
                \expandafter\eeaa\the\fontdimen1\scriptfont1 \ht0 }
\def\accentadjscriptscript#1{\setbox0\hbox{$#1$}\kern   
                \expandafter\eeaa\the\fontdimen1\scriptscriptfont1 \ht0 }
\def\accentadjtextback#1{\setbox0\hbox{$#1$}\kern       
                -\expandafter\eeaa\the\fontdimen1\textfont1 \ht0 }
\def\accentadjscriptback#1{\setbox0\hbox{$#1$}\kern     
                -\expandafter\eeaa\the\fontdimen1\scriptfont1 \ht0 }
\def\accentadjscriptscriptback#1{\setbox0\hbox{$#1$}\kern 
                -\expandafter\eeaa\the\fontdimen1\scriptscriptfont1 \ht0 }

\newcommand{\dmu}{d\mu}

\newcommand{\Om}{\Omega}

\renewcommand{\phi}{\varphi}
\newcommand{\p}{{$p\mspace{1mu}$}}
\newcommand{\R}{\mathbb{R}}

\newcommand{\limminus}{{\mathchoice{\raise.17ex\hbox{$\scriptstyle -$}}
                {\raise.17ex\hbox{$\scriptstyle -$}}
                {\raise.1ex\hbox{$\scriptscriptstyle -$}}
                {\scriptscriptstyle -}}}
\newcommand{\limplus}{{\mathchoice{\raise.17ex\hbox{$\scriptstyle +$}}
                {\raise.17ex\hbox{$\scriptstyle +$}}
                {\raise.1ex\hbox{$\scriptscriptstyle +$}}
                {\scriptscriptstyle +}}}
\newcommand{\limpm}{{\mathchoice{\raise.17ex\hbox{$\scriptstyle \pm$}}
                {\raise.17ex\hbox{$\scriptstyle \pm$}}
                {\raise.16ex\hbox{$\scriptscriptstyle \pm$}}
                {\scriptscriptstyle \pm}}}
\newcommand{\Np}{N^{1,p}}

\newcommand{\g}{\gamma}

\makeatletter
\newcommand{\setcurrentlabel}[1]{\def\@currentlabel{#1}}
\makeatother
\numberwithin{equation}{section}
\newcommand{\phii}[0]{\varphi}

\newcommand{\avint}{\vint}

\begin{document}

\title[Local behavior of \p-harmonic Green's functions]
{Local behavior of \p-harmonic Green's functions in metric spaces}

\author{Donatella Danielli}
\address[Donatella Danielli]{Department of Mathematics \\
Purdue University \\ West Lafayette, IN 47907, USA} 
\email{danielli@math.purdue.edu}

\author{Nicola Garofalo}
\address[Nicola Garofalo]{Department of Mathematics\\Purdue University \\
West Lafayette, IN 47907, USA} 
\email{garofalo@math.purdue.edu}
\thanks{Second author supported in part by NSF Grant DMS-0701001}

\author{Niko Marola}
\address[Niko Marola]{Department of Mathematics and Systems Analysis \\
Helsinki University of Technology \\
P.O. Box 1100 FI-02015 TKK \\
Finland} \email{niko.marola@tkk.fi}
\thanks{Third author supported by the Academy of Finland and Emil
Aaltosen s\"a\"ati\"o}

\date{}

\subjclass[2000]{31C45, 35J60}

\keywords{Capacity, doubling measure, Green function, metric space, 
Newtonian space, \p-Dirichlet integral, \p-harmonic, \p-Laplace equation,
minimizer, Poincar\'e inequality, singular function, Sobolev space}

\maketitle

\emph{This paper is dedicated to the memory of Professor Juha Heinonen}

\bigskip

\begin{abstract}
We describe the behavior of \p-harmonic Green's functions 
near a singularity in metric measure spaces equipped with a 
doubling measure and supporting a Poincar\'e inequality.
\end{abstract}

\section{Introduction}

Holopainen and Shanmugalingam \cite{HoSha} constructed in the metric 
measure space setting a \p-harmonic Green's function, 
called a singular function there,  having most of the 
characteristics of the Green function which is the fundamental solution of
the Laplace operator. 

In this paper we study the following question related to the 
local bahavior of a \p-harmonic Green's function
on locally doubling metric measure
space $X$ supporting a local $(1,p)$-Poincar\`e inequality: Given
a relatively compact domain $\Om \subset X$, $x \in \Om$, and a 
\p-harmonic Green's function $G$ with a singularity at $x$, then can we 
describe the behavior of $G$ near $x$? 

Capacitary estimates for metric rings 
play an important role in the study of the asymptotic
behavior. Following the ideas in the works of Serrin \cite{Serrin1}, 
\cite{Serrin2}, (see also \cite{LSW}) 
such estimates were used in Capogna et al.~\cite{CaDaGa}
to establish the local behavior of singular solutions to a large class of 
nonlinear subelliptic equations which arise in the Carnot--Carath\'eodory geometry. Sharp capacitary estimates for metric rings with unrelated radii
were established in the metric measure space setting
in \cite{GaMa}.

Here, we confine ourselves to mention that a fundamental example of
the spaces included in this paper is obtained by endowing a
connected Riemannian manifold $M$ with the Carath\'eodory metric
$d$ associated with a given subbundle of the tangent bundle,
see \cite{Ca}. If such subbundle generates the tangent space at
every point, then thanks to the theorem of Chow \cite{Chow} and
Rashevsky \cite{Ra} $(M,d)$ is a metric space. Such metric spaces
are known as sub-Riemannian or Carnot-Carath\'eodory (CC) spaces. By
the fundamental works of Rothschild and Stein \cite{RS}, Nagel,
Stein and Wainger \cite{NSW}, and of Jerison \cite{J}, every CC
space is locally doubling, and it locally satisfies a
$(p,p)$-Poincar\`e inequality for any $1\leq p<\infty$. Another
basic example is provided by a Riemannian manifold $(M^n,g)$ with
nonnegative Ricci tensor. In such case thanks to the Bishop
comparison theorem the doubling condition holds globally, see e.g.
\cite{Ch}, whereas the $(1,1)$-Poincar\`e inequality was proved by
Buser \cite{Bu}. An interesting example to which our results apply
and that does not fall in any of the two previously mentioned
categories is the space of two infinite closed cones 
$X=\{(x_1,\ldots, x_n)\in\R^n:\ x_1^2+\ldots +x_{n-1}^2\leq x_n^2\}$ 
equipped with the Euclidean metric of $\R^n$ and with the Lebesgue measure.
This space is Ahlfors regular, and 
it is shown in Haj\l asz--Koskela~\cite[Example 4.2]{HaKo} that a
$(1,p)$-Poincar\'e inequality holds in $X$ if and only if $p>n$. 
Another example is obtained by gluing two copies of closed $n$-balls 
$\{x\in \R^n:\ |x|\leq 1\}$, $n\geq 3$, along a line segment. In this
way one obtains
an Ahlfors regular space that supports a $(1,p)$-Poincare inequality
for $p>n-1$. 
A thorough overview of analysis on metric spaces 
can be found in Heinonen~\cite{heinonen}. One should also consult
Semmes~\cite{Semmes} and David and Semmes~\cite{DaSem}.

Our main result in this paper is a quantative description of
the local behavior of a \p-harmonic Green's function defined in 
Holopainen--Shanmugalingam \cite{HoSha}. 
We shall prove that a Green's function $G$ with a singularity 
at $x_0$ in a relatively compact domain satisfies the asymptotic behavior
\[
G(x) \approx \biggl(\frac{d(x,x_0)^{p}}{\mu(B(x_0,d(x,x_0)))}\biggr)^{1/(p-1)},
\]
where $x$ is uniformly close to $x_0$. Our approach uses 
upper gradients \'a la Heinonen and Koskela~\cite{HeKo}, and 
\p-harmonic functions that can be characterized in terms of 
 \p-energy minimizers among functions with the same boundary
values in relatively compact subsets. Following 
\cite{HoSha} we adopt a definition for Green's functions 
that uses an equation for \p-capacities of level sets.

We want to stress the fact that even in Carnot groups of homogeneous 
dimension $Q$ it is not known whether such \p-harmonic Green's function 
is unique when $1<p<Q$. However, in the conformal case, i.e. when $p=Q$,
the uniqueness for Green's function for the $Q$-Laplace equation
in Carnot groups was settled by Balogh et al. in \cite{BHT}.

The paper is organized as follows. The second section gathers together
the relevant background such as the definition of doubling measures, upper
gradients, Poincar\'e inequality, Newton--Sobolev spaces, and capacity.
In Section 3 we recall sharp capacitary estimates for metric
rings with unrelated radii proved in Garofalo--Marola~\cite{GaMa}. 
In Section 4 we give the definition of Green's functions.
We establish the local behavior of Green's functions in Section 5,
and we also prove a result on the local integrability of 
Green's functions. Section 6 closes the paper with a result on
the local behavior of Cheeger singular functions. In this section
our approach uses Cheeger gradients (see Cheeger~\cite{Cheeger}) 
emerging from a differentiable structure that the ambient metric 
space admits. In particular, \p-harmonic functions can thus be 
characterized in terms of a weak formulation of the \p-Laplace 
equation.

\subsection*{Acknowledgements}
The authors would like to thank Nageswari Shanmugalingam for valuable 
comments on the manuscript and her interest in the paper.

The paper was completed while the third author was visiting Purdue 
University in 2007--2008. He thanks the Department of Mathematics 
for the hospitality and several of its faculty for fruitful
discussions.

\section{Preliminaries}
\label{prelim}

We begin by stating the main assumptions we make on the metric space $X$
and the measure $\mu$.

\subsection{General Assumptions} \label{assumptions}

Throughout the paper $X=(X,d,\mu)$ is a locally compact metric space
endowed with a metric $d$ and a positive Borel regular measure $\mu$.
 We assume that for every compact set $K\subset X$ there exist 
constants $C_K \geq 1$, $R_K>0$ and $\tau_K \ge 1$,
such that for any $x\in K$ and every $0<r\leq R_K$,
$0<\mu(B) <\infty$, where $B:=B(x,r):=\{y\in X:\ d(y,x)<r\}$, 
and, in particular, one has:
\begin{itemize}
\item[(i)] the closed balls $\overline B(x,r)=\{y\in X:d(y,x)\leq r\}$ are compact;
\item[(ii)] (local doubling condition) $\mu(B(x,2r)) \le C_K \mu(B(x,r))$;
\item[(iii)] (local weak $(1,p_0)$-Poincar\'e inequality) there exists $1<p_0<\infty$ such that for all measurable functions $u$ on $X$ and all upper gradients
$g_u$ (see Section 2.3) of $u$
\begin{equation*} \label{PI-ineq}
        \vint_{B(x,r)} |u-u_{B(x,r)}| \,\dmu
        \le C_K r \Big( \vint_{B(x,\tau_K r)} g_u^{p_0} \,\dmu \Big)^{1/p_0},
\end{equation*}
where 
$ u_{B(x,r)} :=\vint_{B(x,r)}u \, d\mu :=\int_{B(x,r)} u\, d\mu/\mu(B(x,r))$.
\item[(iv)] (X is $\LLC$, i.e. linearly locally connected) 
there exists a constant $\alpha \geq 1$ such that for all balls 
$B(x,r)\subset X$, $0<r\leq R_K$, each pair of distinct points 
in the annulus 
$B(x,2r)\setminus\overline{B}(x,r)$
can be connected by a rectifiable path in the annulus 
$B(x,2\alpha r)\setminus\overline{B}(x,r/\alpha)$. 
\end{itemize}

Hereafter, the constants $C_K, R_K$ and $\tau_K$ will be referred to
as the \emph{local parameters} of $K$. We also say that a constant $C$ depends
on the local doubling constant of $K$ if $C$ depends on $C_K$.

The above assumptions encompass, e.g., all Riemannian manifolds with
Ric $\geq 0$, but they also include all Carnot--Carath\'eodory
spaces, and therefore, in particular, all Carnot groups. For a
detailed discussion of these facts  we refer the reader to the paper
by  Garofalo--Nhieu~\cite{GaNhi}. In the case of
Carnot--Carath\'eodory spaces, recall that if the Lie algebra
generating vector fields grow at infinity faster than linearly, then
the compactness of metric balls of large radii may fail in general.
Consider for instance in $\R$ the smooth vector field of H\"ormander
type $X_1=(1+x^2)\frac{d}{dx}$. Some direct calculations prove that
the distance relative to $X_1$ is given by
$d(x,y)=|\arctan(x)-\arctan(y)|$, and therefore, if $r\geq \pi/2$,
we have $B(0,r)=\R$.

\subsection{Local doubling property}
We note that assumption (ii) implies that for every compact set
$K\subset X$ with local parameters $C_K$ and $R_K$, for any $x\in K$
and every $0<r\leq R_K$, one has for $1\leq \lambda \leq R_K/r$,
\begin{equation}\label{dc}
\mu(B(x,\lambda r)) \leq C\lambda^Q\mu(B(x,r)), \end{equation}
 where
$Q=\log_2C_K$, and the constant $C$ depends only on the local doubling
constant $C_K$. The exponent $Q$ serves as a local dimension of the
doubling measure $\mu$ restricted to the compact set $K$.

For $x\in X$ we define the \emph{pointwise dimension} $Q(x)$ by
\begin{align*}
Q(x) = \sup\{& q > 0:\ \exists C>0\ \textrm{ such that } \\
&\lambda^q\mu(B(x,r)) \leq C\mu(B(x,\lambda r)), \\
& \textrm{ for all } 1\leq\lambda<\diam X \textrm{ and } 0<r<\infty\}.
\end{align*}

The inequality \eqref{dc} readily implies that $Q(x) \leq Q$ for
every $x\in K$. Moreover, it follows that
\begin{equation} \label{lowerbound}
\lambda^{Q(x)}\mu(B(x,r)) \leq C\mu(B(x,\lambda r))
\end{equation}
for any $x \in K$, $0<r\leq R_K$ and $1\leq \lambda \leq R_K/r$,
and the constant $C$ depends on the local doubling constant $C_K$.
Furthermore, for all $0<r\leq R_K$ and $x\in K$
\begin{equation} \label{bounds}
C_1r^Q \leq \frac{\mu(B(x,r))}{\mu(B(x,R_K))} \leq C_2r^{Q(x)},
\end{equation}
where $C_1=C(K, C_K)$ and $C_2= C(x,K,C_K)$.

For more on doubling measures, see, e.g. Heinonen~\cite{heinonen}
and the references therein.

\subsection{Upper gradients}
A nonnegative Borel function $g$ on $X$ is an \emph{upper gradient}
of an extended real valued function $f$
on $X$ if for all rectifiable paths $\gamma$
joining points $x$ and $y$ in $X$ we have
\begin{equation} \label{ug-cond}
|f(x)-f(y)|\le \int_\gamma g\,ds.
\end{equation}
whenever both $f(x)$ and $f(y)$ are finite, and
$\int_\g g\, ds=\infty $ otherwise. See Cheeger~\cite{Cheeger}, 
Shanmugalingam~\cite{Sh-rev}, and Heinonen--Koskela~\cite{HeKo} 
for a discussion on upper gradients.

If $g$ is a nonnegative measurable function on $X$
and if (\ref{ug-cond}) holds for \p-almost every path,
then $g$ is a \emph{weak upper gradient} of~$f$.
By saying that (\ref{ug-cond}) holds for \p-almost every path
we mean that it fails only for a path family with zero \p-modulus 
(see, for example, \cite{Sh-rev}).

If $f$ has an upper gradient in $L^p(X)$, then it has a \emph{minimal
weak upper gradient} $g_f \in L^p(X)$ in the sense  that for
every \p-weak upper gradient $g \in L^p(X)$ of $f$,  $g_f \le g$
$\mu$-almost everywhere (a.e.), see Corollary~3.7 in
Shanmugalingam~\cite{Sh-harm}.  The minimal weak upper
gradient can be obtained by the formula
\[
    g_f(x)
     := \inf_g \limsup_{r\to0\limplus}\avint_{B(x,r)} g\,\dmu,
\]
where the infimum is taken over all upper gradients $g \in L^p(X)$ of $f$,
see Lemma~2.3 in Bj\"orn~\cite{Bj}.

\medskip

\subsection{Capacity}
Let $\Om \subset X$ be open and $K \subset \Om$ compact. 
The \emph{relative \p-capacity} of $K$ with respect to $\Om$ is the number 
\begin{equation*} 
  \Capc_p (K,\Om) =\inf\int_\Om g_u^p\,d\mu,
\end{equation*}
where the infimum is taken over all functions $u \in \Np(X)$ such that
$u=1$ on $K$ and $u=0$ on $X\setminus\Om$. If such functions do
not exist, we set $\Capc_p (K,\Om)=\infty$. When $\Om=X$ we simply
write $\Capc_p(K)$.

Observe that the infimum above could be taken over all functions
$u \in \Lip_0(\Om)=\{f \in \Lip(X):\ f=0 \textrm{ on } X\setminus\Om\}$ 
such that $u=1$ on $K$. In addition, the relative \p-capacity is a 
Choquet capacity and consequently for all Borel sets $E$ we have
\[
\Capc_p (E,\Om) = \sup\{\Capc_p(K):\ K\subset E,\ K\textrm{ compact}\}.
\]
For other properties as well as equivalent definitions of
the capacity we refer to Kilpel\"ainen et al.~\cite{KiKiMa},
Kinnunen--Martio~\cite{KiMa96, KiMaNov}, and 
Kallunki--Shanmugalingam~\cite{KaSh}. 

Finally, we say that a property
holds \emph{\p-quasieverywhere} if the set of points for which the property
does not hold is of zero capacity.

\subsection{Newtonian spaces}

We define Sobolev spaces on the metric space
following Shanmugalingam~\cite{Sh-rev}. Let $\Om \subseteq X$ be nonempty and
open. Whenever $u\in L^p(\Om)$, let
$$
        \|u\|_{\Np(\Om)} = \biggl( \int_\Om |u|^p \, \dmu
                + \inf_g  \int_\Om g^p \, \dmu \biggr)^{1/p},
$$
where the infimum is taken over all weak upper gradients of $u$.
The \emph{Newtonian space} on $\Om$ is the quotient space
$$
        \Np (\Om) = \{u: \|u\|_{\Np(\Om)} <\infty \}/{\sim},
$$
where  $u \sim v$ if and only if $\|u-v\|_{\Np(\Om)}=0$.
The Newtonian space is a Banach space and a lattice, moreover
Lipschitz functions are dense, see \cite{Sh-rev} and Bj\"orn
et al.~\cite{BBS}.

To be able to compare the boundary values of Newtonian functions
we need a Newtonian space with zero boundary values.
Let $E$ be a measurable subset of $X$. The \emph{Newtonian space with zero boundary values} is the space
\[
\Np_0(E)=\{u|_{E} : u \in \Np(X) \text{ and }
        u=0 \text{ on } X \setm E\}.
\]
The space $\Np_0(E)$ equipped with the norm inherited from $\Np(X)$
is a Banach space, see Theorem~4.4 in Shanmugalingam~\cite{Sh-harm}.

We say that $u$ belongs to the \emph{local Newtonian space}
$\Np\loc(\Omega)$ if $u\in \Np(\Om')$ for every open
$\Om'\Subset\Omega$ (or equivalently that $u\in \Np(E)$ for every measurable $E\Subset\Omega$).

We will also need an inequality for Newtonian functions with zero
boundary values. If $f\in\Np_0(B(x,r))$, then there exists a constant
$C>0$ only depending on $p$, the local doubling constant, 
and the constants in the weak Poincar\'e inequality, such that
\begin{equation} \label{eq:SoboPI}
\biggl(\avint_{B(x,r)}|f|^{p}\, d\mu\biggr)^{1/p} \leq 
Cr\biggl(\avint_{B(x,r)}g_f^p\,d\mu\biggr)^{1/p}
\end{equation}
for every ball $B(x,r)$ with $r\leq \frac1{3}\diam X$.
For this result we refer to Kinnunen and Shanmugalingam~\cite{KiSh1}.

\subsection{Differentiable structure}

Cheeger~\cite{Cheeger} demonstrated that metric measure spaces that satisfy
assumptions (ii) and (iii) admit a differentiable structure with which Lipschitz functions can be differentiated almost everywhere. This differentiable 
structure gives rise to an alternative definition of a Sobolev space 
over the given metric measure space than defined above. However, 
assuming (ii) and (iii) these definitions lead to the same space, 
see Shanmugalingam~\cite[Theorem 4.10]{Sh-rev}.
Thanks to a deep theorem by Cheeger the corresponding Sobolev 
space is reflexive, see 
\cite[Theorem 4.48]{Cheeger}. 

The differentiable structure gives the notion of partial 
derivatives in the following theorem, see 
Cheeger~\cite[Theorem 4.38]{Cheeger},
and it is compatible with the notion of an upper gradient.

\medskip

\begin{thm}[Cheeger]
Let $X$ be a metric measure space equipped with a doubling Borel
regular measure $\mu$. Assume that $X$ admits a weak $(1,p_0)$-Poincar\'e
inequality for some $1<p_0<\infty$. Then there exists measurable sets 
$U_\alpha$ with positive measure such that
\[
\mu(X\setminus\bigcup_\alpha U_\alpha) =0,
\] 
and Lipschitz ``coordinate charts'' 
\[
\mathcal{X}^\alpha = (X_1^\alpha,\ldots,X_{k(\alpha)}^\alpha):X\to\R^{k(\alpha)}
\]
such that for each $\alpha$ functions $X_1^\alpha,\ldots,X_{k(\alpha)}^\alpha$
are linearly independent on $U_\alpha$ and 
\[
1\leq k(\alpha) \leq N, 
\]
where $N$ is a constant depending only on the doubling constant of $\mu$ and
the constants in the Poincar\'e inequality. Moreover, if $f:X\to\R$ is
Lipschitz, then there exist unique (up to a set of measure zero) 
bounded vector-valued functions
$d^\alpha f:U_\alpha \to \R^{k(\alpha)}$ such that
\[
\lim_{r\to 0\limplus}\sup_{x\in B(x_0,r)}
\frac{|f(x)-f(x_0)-d^\alpha f\cdot(\mathcal{X}^\alpha(x)-\mathcal{X}^\alpha(x_0))|}{r} = 0
\]
for $\mu$-a.e. $x_0 \in U_\alpha$.
\end{thm}

\medskip

We can assume that the sets $U_\alpha$ are pairwise disjoint, and
extend $d^\alpha f$ by zero outside $U_\alpha$. Regard $d^\alpha f$ as
vectors in $\R^N$ and let $Df := \sum_\alpha d^\alpha f$. 
By Shanmugalingam~\cite[Theorem 4.10]{Sh-rev} and 
\cite[Theorem 4.47]{Cheeger}, the Newtonian space $N^{1,p_0}(X)$ is
equal to the closure in the $N^{1,p_0}$-norm of the collection of (locally)
Lipschitz functions on $X$, then the derivation operator $D$ can
be extended to all of $N^{1,p_0}(X)$ so that there exists a constant 
$C>0$ such that
\[
C^{-1}|Df(x)| \leq g_f(x) \leq C|Df(x)|
\]
for all $f\in N^{1,p_0}(X)$ and $\mu$-a.e. $x\in X$.
Here the norms $|\cdot|$ can be chosen to be inner product norms.
The differential mapping $Df$ satisfies
the product and chain rules: if $f$ is a bounded Lipschitz function
on $X$, $u\in N^{1,p_0}(X)$, and $h:\R\to\R$ is continuously differentiable
with bounded derivative, then $uf$ and $h\circ u$ both belong to 
$ N^{1,p_0}(X)$ and
\begin{align*}
D(uf) & = uDf + fDu; \\
D(h\circ u) = (h\circ u)'Du.
\end{align*}
See the discussion in Cheeger~\cite{Cheeger} and Keith~\cite{Keith}.

\subsection{\p-harmonic functions}
Let $\Om\subset X$ be a domain. A
function $u \in \Np\loc(\Om)\cap C(\Om)$ is \emph{\p-harmonic} in $\Om$ if
for all relatively compact sets $\Om'\subset\Om$ and for all $\phii \in \Np_0(\Om')$,
\[
\int_{\Om'}g_u^p\,d\mu \leq \int_{\Om'}g_{u+\phii}^p\,d\mu.
\]
It is known that nonnegative \p-harmonic functions satisfy
Harnack's inequality and the strong maximum principle, there
are no non-constant nonnegative \p-harmonic functions on all of $X$,
and \p-harmonic functions have locally H\"older continuous
representatives. See \cite{KiSh1}.

As a consequence of the $\LLC$ property of $X$ a nonnegative
\p-harmonic function on an annulus $B(y,Cr)\setminus B(y,r/C)$
satisfies Harnack's inequality on the sphere
$S(y,r)=\{x\in X:\ d(x,y) = r\}$ for sufficiently small $r$, see
Bj\"orn et al. \cite[Lemma 5.3]{BMcSha}.

We also say that a function  $u \in \Np\loc(\Om)\cap C(\Om)$ is 
\emph{Cheeger \p-harmonic} in $\Om$ if in the above definition upper
gradients $g_u$ and $g_{u+\phii}$ are replaced by $|Du|$ and
$|D(u+\phii)|$, respectively. 
Note that by a result in Cheeger~\cite{Cheeger}, the Cheeger
\p-harmonic functions 
are \p-quasiminimizers in the sense of, e.g.,
Kinnunen--Shanmugalingam~\cite{KiSh1}. Moreover, the Cheeger \p-harmonic 
functions can be characterized in terms of a weak 
formulation of the \p-Laplace equation: $u$ is Cheeger 
\p-harmonic if and only if
\[
\int_{\Om'}|Du|^{p-2}Du\cdot D\phii\, d\mu = 0
\]
for all $\Om'$ and $\phii$ as in the above definition.

\section{Capacitary estimates}

The aim of this section is to recall sharp capacity estimates for 
metric rings with unrelated radii proved in \cite{GaMa}. 
We emphasize an interesting feature of Theorems \ref{thm:below} and 
\ref{thm:above} that cannot be observed in the setting of, for example, 
Carnot gouprs. That is the dependence of the estimates on the center
of the ring. This is a consequence of the fact that in the general
setting $Q(x_0) \neq Q$ where $x_0 \in X$, see Section~\ref{prelim}.
The results in this section will play an important role in the subsequent 
developments.

\medskip
 
\begin{thm}(Estimates from below) \label{thm:below}
Let $\Om \subset X$ be a bounded open set, $x_0 \in \Om$, and $Q(x_0)$ be
the pointwise dimension at $x_0$. Then there exists $R_0(\Om)>0$
such that for any $0<r<R\leq R_0(\Om)$
we have
\begin{align*}
& \Capc_{p_0}(\overline{B}(x_0,r),B(x_0,R)) \geq \\
& \left\{\begin{array}{ll}
C_1(1-\frac{r}{R})^{p_0(p_0-1)}\frac{\mu(B(x_0,r))}{r^{p_0}},\ \textrm{ if }\ 1<p_0<Q(x_0), \\
C_2(1-\frac{r}{R})^{Q(x_0)(Q(x_0)-1)}\biggl(\log\frac{R}{r}\biggr)^{1-Q(x_0)},
\ \textrm{ if }\ p_0=Q(x_0), \\
C_3(1-\frac{r}{R})^{p_0(p_0-1)}\biggl|(2R)^{\frac{p_0-Q(x_0)}{p_0-1}}-r^{\frac{p_0-Q(x_0)}{p_0-1}}\biggr|^{1-p_0},\ \textrm{ if }\ p_0>Q(x_0),
\end{array} \right.
\end{align*}
where
\begin{align*}
C_1 & = C\biggl(1-\frac1{2^{\frac{Q(x_0)-p_0}{p_0-1}}}\biggr)^{p_0-1}, \\
C_2 & =C\frac{\mu(B(x_0,r))}{r^{Q(x_0)}},\\
C_3 & =C\frac{\mu(B(x_0,r))}{r^{Q(x_0)}}\biggl(2^{\frac{p_0-Q(x_0)}{p_0-1}}-1\biggr)^{p_0-1},
\end{align*}
with $C > 0$ depending only on $p_0$ and the local doubling constant 
of $\Om$.
\end{thm}

\medskip

\begin{remark}
Observe that if $X$ supports the weak $(1,1)$-Poincar\'e inequality,
i.e. $p_0 =1$, these estimates reduce to the capacitary estimates, e.g., 
in Capogna et al.~\cite[Theorem 4.1]{CaDaGa}.
\end{remark}

\medskip

\begin{thm}(Estimates from above) \label{thm:above}
Let $\Om$, $x_0$, and $Q(x_0)$ be as in Theorem~\ref{thm:below}.
Then there exists $R_0(\Om)>0$ such that for any $0<r<R\leq R_0(\Om)$ we have
\begin{align*}
& \Capc_{p_0}(\overline{B}(x_0,r),B(x_0,R)) \\
& \leq \left\{\begin{array}{ll}
C_4\frac{\mu(B(x_0,r))}{r^{p_0}}, & \textrm{ if }\,
1<p_0<Q(x_0), \\
C_5\biggl(\log\frac{R}{r}\biggr)^{1-Q(x_0)}, & \textrm{ if }\, p_0=Q(x_0), \\
C_6\left|(2R)^{\frac{p_0-Q(x_0)}{p_0-1}}-r^{\frac{p_0-Q(x_0)}{p_0-1}}\right|^{1-p_0}, &
\textrm{ if }\, p_0>Q(x_0),
\end{array} \right.
\end{align*}
where $C_4$ is a positive constant depending only
on  $p_0$ and the local doubling constant of $\Om$, whereas
\[
C_5 = C\frac{\mu(B(x_0,r))}{r^{Q(x_0)}},
\]
where $C$ is a positive constants depending only
on  $p_0$ and the local doubling constant of $\Om$, and
\[
C_6 =C\biggl(2^{\frac{p_0-Q(x_0)}{p_0-1}}-1\biggr)^{-1},
\]
with $C > 0$ depending on $p_0$, the local parameters of $\Om$, and 
$\mu(B(x_0,R_0))$. 
\end{thm}

\medskip

We have the following immediate corollary.

\medskip

\begin{cor}
If $1< p_0 \leq Q(x_0)$, then we have
\begin{equation*} \label{eq:singleton}
\Capc_{p_0}(\{x_0\},\Om) = 0.
\end{equation*}
\end{cor}

\section{Green's functions}

We define a Green's function on metric spaces following Holopainen and
Shanmugalingam~\cite{HoSha}. Note that 
Holopainen and Shanmugalingam referred to this function class 
as singular functions. We consider here a definition that
uses an equation for \p-capacities of level sets. Green's function
on a Riemannian manifold satisfies this equation, see 
Holopainen~\cite{Holopainen}.

\medskip

\begin{deff} \label{def:G}
Given $1<p_0\leq Q(x_0)$, let $\Om \subset X$ be a
relatively compact domain, and
$x_0 \in \Om$. An extended real-valued function $G = G(\cdot,x_0)$
on $\Om$ is said to be a \emph{Green's function with singularity at $x_0$}
if the following criteria are satisfied:

\medskip

\begin{enumerate}

\item[1.]
$G$ is \p$_0$-harmonic and positive in $\Om\setminus\{x_0\}$,

\item[2.]
$G|_{X\setminus\Om}=0$ \p-quasieverywhere and 
$G\in N^{1,p_0}_{\loc}(X\setminus B(x_0,r))$ for all $r>0$,

\item[3.]
$x_0$ is a singularity, i.e.,
\[
\lim_{x\to x_0}G(x) = \infty.
\]

\item[4.]
whenever $0\leq\alpha < \beta$,
\[
C_1(\beta-\alpha)^{1-p_0}\leq \Capc_{p_0}(\Om^\beta,\Om_\alpha)) \leq
C_2(\beta-\alpha)^{1-p_0},
\]
where $\Om^\beta = \{x\in\Om:\ G(x)\geq \beta\}$,
$\Om_\alpha = \{x\in\Om:\ G(x) > \alpha\}$, and $C_1,\,C_2 >0$ are
constants depending only on $p_0$.
\end{enumerate}
\end{deff}

\medskip

\begin{remark}(Existence)
The existence of Green's functions in the 
$Q$-regular metric space setting was first proved by Holopainen and Shanmugalingam in \cite{HoSha}. Being
a $Q$-regular metric measure space means that the measure $\mu$
satisfies, for all balls $B(x,r)$ a double inequality
\[
C^{-1}r^Q\leq \mu(B(x,r))\leq Cr^Q
\]
with a fixed constant $Q$. There are, however, many instances 
where the $Q$-regularity condition is not satisfied. 
For example, systems of vector fields of H\"ormander type are, 
in general, not $Q$-regular for any $Q>0$. 

In \cite{GaMa} the $Q$-regularity assumption was removed and the
existence of this function class was proved in more general setting. 
For the proof of the existence,
we refer to \cite[Theorem 3.4]{HoSha}, see also remarks in \cite{GaMa}. 
\end{remark}

\medskip

\begin{remark}(Uniqueness)
It is not known whether a Green's function is unique in the metric space 
setting even in the case of Cheeger \p-harmonic functions. 
Indeed, the uniqueness of Green's functions 
is not settled in Carnot groups when $1<p_0<Q$, where $Q$ is the 
homogeneous dimension attached to the non-isotropic dilations. However,
Green's function is known to be unique when $p_0=Q$, 
see Balogh et al.~\cite{BHT}. 
\end{remark}

\section{Local behavior of \p-harmonic Green's functions}

We begin by recalling that if $K\subset\Om$ is closed, $u\in N^{1,p_0}(X)$ is
a \emph{\p$_0$-potential} of $K$ (with respect to $\Om$) if 
\begin{itemize}

\item[(i)] $u$ is \p$_0$-harmonic on $\Om\setminus K$; 

\item[(ii)] $u=1$ on $K$ and $u=0$ in $X\setminus\Om$. 

\end{itemize}

By Lemma 3.3 in Holopainen--Shanmugalingam~\cite{HoSha} 
\p$_0$-potentials always exist if $\Capc_{p_0}(K,\Om)<\infty$.

\medskip

From know on, we set
\[
m(r) = m_G(x_0,r) = \min_{\partial B(x_0,r)}G, \quad M(r)=M_G(x_0,r) 
= \max_{\partial B(x_0,r)}G,
\]
where $G$ is a Green's function with singularity at $x_0$. We can now 
state the following growth estimates for a Green's
function near a singularity. 
In what follows, $R_0(\Om) >0$ is the constant from
theorems \ref{thm:below} and \ref{thm:above}.

\medskip

\begin{thm} \label{thm:blowupI}
Let $\Om$ be a relatively compact domain in $X$, $x_0 \in \Om$, and
$1<p\leq Q(x_0)$. If
$G$ is a Green's function with singularity at $x_0$ and given
$0<R\leq R_0(\Om)$ for which $\overline{B}(x_0,R) \subset \Om$, then
for every $0<r<R$ we have
\begin{equation*}
m_G(x_0,r) \leq C_1\biggl(\frac1{\Capc_{p_0}(\overline{B}(x_0,r),B(x_0,R))}
\biggr)^{1/(p_0-1)} +M_G(x_0,R).
\end{equation*}
Suppose $r_0 \in (0,R)$ is such that $m_G(x_0,r_0)\geq M_G(x_0,R)$, then
for every $0<r<r_0$ we have
\begin{equation*}
M_G(x_0,r) \geq C_2\biggl(\frac1{\Capc_{p_0}(\overline{B}(x_0,r),B(x_0,r_0))}
\biggr)^{1/(p_0-1)} + M_G(x_0,R),
\end{equation*}
where the constants $C_1$ and $C_2$ both depend only on $p_0$.
\end{thm}

\medskip

\begin{proof}
Consider a radius $R>0$ such that $\overline{B}(x_0,R)\subset\Om$. 
Since $G(x) \to \infty$ when $x$ tends to $x_0$, 
the maximum principle implies that
\begin{equation} \label{eq:minest}
m(r) \geq m(\rho), \quad 0<r<\rho<R.
\end{equation}
Define $w=G-M(R)$, and hence $w\leq 0$ on $\partial B(x_0,R)$. 
Observe that the first inequality in the theorem obviously holds true
if $m(r) \leq M(R)$, thus, we might as well assume that
\begin{equation} \label{assumption1}
m(r) > M(R),
\end{equation}
and consider the function $v$ 
in the annulus $B(x_0,R)\setminus\overline{B}(x_0,r)$ defined by
\begin{equation*}
v = \left\{\begin{array}{ll}
0, & \textrm{ if }\ G\leq M(R), \\
w, & \textrm{ if }\ M(R) < G < m(r), \\
m_w(r), & \textrm{ if }\ G \geq m(r).
\end{array} \right.
\end{equation*}
If we extend $v$ by letting $v=m_w(r)$ on $\overline{B}(x_0,r)$, then 
$v\in N^{1,p_0}_0(B(x_0,R))$. Our assumption \eqref{assumption1} implies that
$m_w(r)=m(r)-M(R) > 0$, so the function
\[
\phii = \frac{v}{m_w(r)},
\]
which equals to $1$ in $\overline{B}(x_0,r)$, is both an admissible 
function for the capacity of $\overline{B}(x_0,r)$ with respect to 
$B(x_0,R)$ and the \p$_0$-potential of the set $\{x\in X:\ \phii(x)\geq 1\}$ 
with respect to the set $\{x\in X:\ \phii(x)>0\}$. Thus one has
\begin{align*}
& \Capc_{p_0}(\overline{B}(x_0,r),B(x_0,R)) \leq \int_{B(x_0,R)}
g_{\phii}^{p_0}\,d\mu \\
& = \Capc_{p_0}(\{x\in X:\ \phii(x)\geq 1\},\{x\in X:\ \phii(x)>0\}) \\
& = \Capc_{p_0}(\{x\in X:\ G(x)\geq m(r)\},\{x\in X:\ G(x)>M(R)\}) \\  
& \leq C_1(m(r)-M(R))^{1-p_0},
\end{align*}
where we used criterion 4 from Definition~\ref{def:G} and the fact that
$\phii \geq 1$ or $\phii >0$ if and only if $G\geq m(r)$ or $G>M(r)$,
respectively. This implies the first claim.

To prove the second iequality of the claim, let $w=G-M(R)$. 
Let $r_0\in(0,R)$ be such that $m(r_0)\geq M(R)$. This implies that
$w\geq 0$ on $\overline{B}(x_0,r_0)$ and also that $M(r)\geq M(R)$, 
for all $0<r<r_0$. Hence, by the maximum principle
we have that
\[
\{x\in \Om:\ G(x)\geq M(r)\} \subset \overline{B}(x_0,r)
\]
and
\[
B(x_0,r_0) \subset \{x\in\Om:\ G(x)>M(R)\}.
\]
Hence it follows that
\begin{align*}
& \Capc_{p_0}(\overline{B}(x_0,r), B(x_0,r_0)) \\
& \geq \Capc_{p_0}(\{x\in \Om:\ G(x)\geq M(r)\}, B(x_0,r_0)) \\
& \geq \Capc_{p_0}(\{x\in \Om:\ G(x)\geq M(r)\},\{x\in\Om:\ G(x)>M(R)\}) \\
& \geq C_2(M(r)-M(R))^{1-p_0}, 
\end{align*}
which implies the second claim and the proof is complete.
\end{proof}

\medskip

\begin{remark}
Theorem 7.1 in Capogna et al.~\cite{CaDaGa} is slightly 
incorrect as the additional term $M(R)$ is missing from
the left-hand and the right-hand side in (ii). However, this does  
not affect the results in that paper since the additional term can be
absorbed when establishing results on the behavior near a
singularity.
\end{remark}

\medskip

We have the following result on the local behavior of a Green's
function near a singularity. 

\medskip

\begin{thm} \label{thm:localbI}
Let $\Om$ be a relatively compact domain in $X$, and $x_0 \in \Om$.
If $G$ is a Green's function with singularity at $x_0$, then 
there exist positive constants $C_1,C_2$ and 
$R_0$ such that for any $0<r<\frac{R_0}{2}$ and $x \in B(x_0,r)$ we have
\begin{align*}
C_1\biggl(\frac{d(x,x_0)^{p_0}}{\mu(B(x_0,d(x,x_0)))}\biggr)^{1/(p_0-1)} &
\leq G(x) \\
& \leq C_2\biggl(\frac{d(x,x_0)^{p_0}}{\mu(B(x_0,d(x,x_0)))}\biggr)^{1/(p_0-1)},
\end{align*}
when $1<p_0<Q(x_0)$, whereas
\[
C_1 \log\left(\frac{R_0}{d(x,x_0)}\right) \leq G(x) 
\leq C_2\log\left(\frac{R_0}{d(x,x_0)}\right),
\]
when $p_0=Q(x_0)$. 
Here the constants $C_1$ and $C_2$ depend on $p_0$, $x_0$, 
and the local parameters of $\Om$, whereas constant $R_0$ depends only 
on $\Om$.
\end{thm}

\medskip

\begin{proof}
Let $R_0=\min\{r_0,R_0(\Om)\}$, where $r_0>0$ is from the second estimate 
in Theorem~\ref{thm:blowupI}. 
The Harnack inequality on a sphere implies that 
there exists a constant $C>0$ such that
\[
M(r) \leq Cm(r).
\]
for every $0< r < R_0$. Let, in particular, $r:=d(x_0,x) < 
\frac{R_0}{2}$. 
From the first estimate in Theorem~\ref{thm:blowup}, 
the maximum principle, and the Harnack
inequality on the sphere, we obtain for any $0< r < \frac{R_0}{2}$
\begin{align*}
G(x) & \leq M(r) \leq Cm(r) \\
& \leq C\Capc_{p_0}(\overline{B}(x_0,r),B(x_0,R_0))^{-1/(p_0-1)}.
\end{align*}
Thanks to Theorem~\ref{thm:below} we have
\begin{align*}
G(x) & \leq C\biggl(1-\frac{r}{R_0}\biggr)^{-p_0}\biggl(\frac{r^{p_0}}{\mu(B(x_0,r))}\biggr)^{1/(p_0-1)} \\
& \leq C\biggl(\frac{r^{p_0}}{\mu(B(x_0,r))}\biggr)^{1/(p_0-1)},
\end{align*}
when $1<p<Q(x_0)$, and
\[
G(x) \leq C\log\left(\frac{R_0}{r}\right),
\]
when $p=Q(x_0)$. This proves the estimate from above. 

To show the estimate from below, observe that the second estimate in 
Theorem~\ref{thm:blowupI}, the maximum principle, and the Harnack
inequality on a sphere imply for $0<r< R_0$
\begin{align*}
G(x) & \geq m(r) \geq C^{-1}M(r) \\
& \geq C\Capc_{p_0}(\overline{B}(x_0,r),B(x_0,R_0))^{-1/(p_0-1)}
\end{align*}
Applying Theorem~\ref{thm:above} we conclude for $1<p_0<Q(x_0)$
\[
G(x) \geq C\biggl(\frac{r^{p_0}}{\mu(B(x_0,r))}\biggr)^{1/(p_0-1)},
\]
and for $p_0=Q(x_0)$ that
\[
G(x) \geq C\log\left(\frac{R_0}{r}\right).
\] 
This completes the proof.
\end{proof}

\medskip

\begin{remark}
Note that if $1<p_0<Q(x_0)$ then due to \eqref{bounds}, 
it readily follows that
\[
C_1d(x,x_0)^{(p_0-Q(x_0))/(p_0-1)} \leq G(x) \leq C_2d(x,x_0)^{(p_0-Q)/(p_0-1)},
\]
when $x \in B(x_0,r)$ with $0<r<\frac{R_0}{2}$. Here the constants
$C_1$ and $C_2$ depend on $p_0,\, x_0$ and the 
local parameters of $\Om$. 
\end{remark}

\medskip

In general Green's function $G\notin L^{p_0}_{\loc}(\Om)$, but as a
corollary of Theorem~\ref{thm:localbI} we have 
the following integrability result near a singularity.

\medskip

\begin{corollary}
Let $1<p_0<Q(x_0)$. Under the assumptions of Theorem~\ref{thm:localb},
one has
\begin{itemize}

\item[(i)]
\[
G \in \bigcap_{0<q<\frac{Q(x_0)(p_0-1)}{Q-p_0}}L^q(B(x_0,r)),
\]

\item[(ii)]
\[
g_G \in \bigcap_{0<q<\frac{Q(x_0)(p_0-1)}{Q-1}}L^q(B(x_0,r)),
\]

\item[(iii)] 
If $p_0 > (Q + Q(x_0) - 1)/Q(x_0)$, then
\[
G \in \bigcap_{1< q<\frac{Q(x_0)(p_0-1)}{Q-1}}N_0^{1,q}(B(x_0,r)).
\]
\end{itemize}
\end{corollary}

\begin{proof}
The proof of (i) is an immediate consequence of the estimate from above 
in Theorem~\ref{thm:localb}. To prove (ii),
we note that since $1<p_0<Q(x_0)\leq Q$, 
\[
q^*:=\frac{Q(x_0)(p_0-1)}{Q-1} < p_0.
\]
Applying H\"older's inequality, the Caccioppoli inequality, see 
Bj\"orn--Marola~\cite[Proposition 7.1]{BjMa},
and again Theorem~\ref{thm:localb}, we find
for $0<q<p$ and for $\sigma \in (0,r)$
\[
\int_{B(x_0,2\sigma)\setminus B(x_0,\sigma)}g_G^q\, d\mu \leq C\sigma^{Q(x_0)-\frac{q(Q-1)}{p_0-1}}.
\]
Note that the exponent $Q(x_0)-\frac{q(Q-1)}{p_0-1}$ is strictly positive,
when $0<q<q^*$ and zero when $q=q^*$. This observation gives us that
\begin{align*}
\int_{B(x_0,r)}g_G^q\, d\mu & = \sum_{i=0}^\infty\int_{B(x_0,2^{-i}r)\setminus B(x_0,2^{-(i+1)}r)}g_G^q\, d\mu \\
& \leq C\mu(B(x_0,r))\sum_{i=0}^\infty(2^{-i}r)^{Q(x_0)-\frac{q(Q-1)}{p_0-1}}< \infty.
\end{align*}
This proves (ii). Finally, (iii) follows from (ii) once we observe that 
the condition $p_0 > (Q+Q(x_0)-1)/Q(x_0)$ is equivalent to $Q(x_0)(p_0-1)/(Q-1) > 1$.
\end{proof}

\section{Cheeger singular functions}

In this section we study Cheeger singular functions, i.e. functions
that satisfy \emph{only} conditions 1, 2 and 3 in Definition~\ref{def:G} and
the notion
of a \p$_0$-harmonic function is replaced by that of a Cheeger 
\p$_0$-harmonic function.

Let $G'$ be a functions that satisfies conditions 1.--3. in 
Definition~\ref{def:G}. We begin by defining $K(G)$ by
\begin{equation} \label{eq:Ku}
K(G') = \int_\Om|DG'|^{p_0-2}DG'\cdot D\phii\, d\mu,
\end{equation}
where $\phii \in N^{1,p_0}_0(\Om)$ is such that $\phii = 1$ in a neighborhood of
$x_0$. If $\phii_i\in N^{1,p_0}_0(\Om), i=1,2,$ and $\phii_i=1$ in a neighborhood of
$x_0$ then $\phii=\phii_1-\phii_2\in N^{1,p_0}_0(\Om\setminus\{x_0\})$. This
gives us
\[
\int_\Om|DG'|^{p_0-2}DG'\cdot D\phii_1\, d\mu = \int_\Om|DG'|^{p_0-2}DG'\cdot
D\phii_2\, d\mu.
\]
Thus $K(G')=K(G',p_0,\Om)$, in particular, $K$ does not depend on
$\phii$. 
Another property of $K(G')$ that will play an important role is that
\[
K(G') > 0,
\]
see \eqref{eq:K=} below. We obtain the following result on the growth
of Cheeger singular functions near a singularity. 

\medskip

\begin{thm} \label{thm:blowup}
Let $\Om$ be a relatively compact domain in $X$, $x_0 \in \Om$, and
$1<p<Q(x_0)$. If
$G'$ is a Cheeger singular function, i.e. $G'$ satisfies conditions
1--3 in Definition~\ref{def:G}, with singularity at $x_0$ and given
$0<R\leq R_0(\Om)$ for which $\overline{B}(x_0,R) \subset \Om$, then
for every $0<r<R$ we have
\begin{equation*}
m_{G'}(x_0,r) \leq \biggl(\frac{K(G')}{\Capc_{p_0}(\overline{B}(x_0,r),B(x_0,R))}\biggr)^{1/(p_0-1)} +M_{G'}(x_0,R).
\end{equation*}
Suppose $r_0 \in (0,R)$ is such that $m_{G'}(x_0,r_0)\geq M_{G'}(x_0,R)$, 
then for every $0<r<r_0$ we have
\begin{multline*}
M_{G'}(x_0,r) \geq C(1-\frac{r}{r_0})^{p_0}\biggl(\frac{K(G')}{\Capc_{p_0}(\overline{B}(x_0,r),B(x_0,r_0))}\biggr)^{1/(p_0-1)} \\ + M_{G'}(x_0,R),
\end{multline*}
where $C=(C_1/C_4)^{1/(p_0-1)} > 0$, and 
the constants $C_1$ and $C_4$ are as in theorems~\ref{thm:below} and 
\ref{thm:above}, respectively.
\end{thm}

\medskip

\begin{proof}
Consider a radius $R>0$ such that $\overline{B}(x_0,R)\subset\Om$. 
Define $w=G'-M(R)$, and hence $w\leq 0$ on $\partial B(x_0,R)$.
Observe also that the first inequality in the theorem obviously 
holds true if $m(r) \leq M(R)$, thus, we might as well assume that
$m(r) > M(R)$.
Let functions $v$ and $\phii=v/m_w(r)$ be defined as in the proof of
Theorem~\ref{thm:blowupI} with $G$ replaced by $G'$. Then $\phii$
can be used in the definition of $K(G')$, see \eqref{eq:Ku}. We have
\begin{align*}
K(G') & = \int_{B(x_0,R)\setminus \overline{B}(x_0,r)}|DG'|^{p_0-2}DG'\cdot 
D\phii\, d\mu \\
& = \frac1{m_w(r)}\int_{B(x_0,R)\setminus \overline{B}(x_0,r)}
|DG'|^{p_0-2}DG'\cdot Dv\, d\mu.
\end{align*}
Observing that $Dv = 0$ whenever $v\neq w$, whereas $Dv=Dw=DG'$ on
the set where $v=w$, we conclude
\begin{equation} \label{eq:K=}
K(G') = \frac1{m_w(r)}\int_{B(x_0,R)}|Dv|^{p_0}\, d\mu = 
m_w^{p_0-1}\int_{B(x_0,R)}|D\phii|^{p_0}\, d\mu.
\end{equation}
Note at this point that \eqref{eq:K=} proves that $K(G')>0$. 
Indeed, if, in fact, $K(G')\leq 0$, the Sobolev--Poincar\'e 
inequality \eqref{eq:SoboPI} implies that 
\[
\int_{B(x_0,R)}|v|^{p_0}\, d\mu \leq CR^{p_0}\int_{B(x_0,R)}
|Dv|^{p_0}\, d\mu \leq 0,
\]
and, moreoover, $v\equiv 0$ in $B(x_0,R)$. This, in turn, 
would contradict the fact that $G'(x) \to \infty$ when
$x$ tends to $x_0$. This shows that $K(G')>0$.

Observing that $\phii=v/m_w(r)$ is an admissible function for the 
capacity of $B(x_0,r)$ with respect to $B(x_0,R)$, we obtain
from \eqref{eq:K=} that
\begin{align} \label{eq:upperbound}
\Capc_{p_0} & (\overline{B}(x_0,r),B(x_0,R)) \leq \int_{B(x_0,R)\setminus 
\overline{B}(x_0,r)}
|D(v/m_w(r))|^{p_0}\, d\mu \\
& \leq \frac1{m_w(r)^{p_0}}\int_{B(x_0,R)\setminus \overline{B}(x_0,r)}|Dv|^{p_0}\, d\mu \leq m_w(r)^{1-p_0}K(G'). \nonumber
\end{align}
This implies the first claim.

To prove the second iequality of the claim, we observe that $w(x) \to 
\infty$, when $x$ tends to $x_0$. As above, $w=G'-M(R)$. 
Also thanks to \eqref{eq:minest} one has that
\[
m_w(r) \geq m_w(\rho), \quad 0<r<\rho<R.
\]
Let $r_0\in(0,R)$ be such that $m(r_0)\geq M(R)$. This implies that
$w\geq 0$ on $\overline{B}(x_0,r_0)$. 
For any $0<r<r_0$ consider the function
$\psi:\R\to\R$ defined by
\[
\psi(t) = \left\{\begin{array}{ll}
1, & \textrm{ in } 0\leq t \leq r, \\
\frac{t^{\frac{p_0-Q(x_0)}{p_0-1}}-r_0^{\frac{p_0-Q(x_0)}{p_0-1}}}{r^{\frac{p_0-Q(x_0)}{p_0-1}}-r_0^{\frac{p_0-Q(x_0)}{p_0-1}}}, & \textrm{ in } 
r\leq t\leq r_0, \\
0, & \textrm{ in } r_0 \leq t \leq R.
\end{array} \right.
\]
Observe that $\psi \in  L^\infty(\R)$, 
$\spt(\psi') \subset [r,r_0]$, and that $\psi' \in L^\infty(\R)$,
thus $\psi$ is a Lipschitz function.
Moreover, $\psi\circ d(x_0,x) \in N^{1,p_0}(B(x_0,R))$. 
As in the proof of Theorem 4.5 in Garofalo--Marola~\cite{GaMa}, we obtain
\[
\int_{B(x_0,R)}|D\psi|^p_0\,d\mu \leq C_4\frac{\mu(B(x_0,r))}{r^{p_0}}.
\]
On the other hand, if we use Theorem~\ref{thm:below}, for the proof 
see \cite{GaMa}, we have
\begin{equation} \label{eq}
\int_{B(x_0,R)}|D\psi|^{p_0}\,d\mu \leq \frac{C_4}{C_1}
(1-\frac{r}{r_0})^{p_0(1-p_0)}\Capc_{p_0}(\overline{B}(x_0,r),B(x_0,r_0)).
\end{equation}
Since $\psi\circ d(x_0,x)$ is an admissible function for $K(G')$, it
follows from \eqref{eq:Ku}, \eqref{eq}, and H\"older's inequality that
\begin{align} \label{eq:upper}
& K(G')^{p_0/(p_0-1)} \\
& \leq \biggl(\int_{B(x_0,R)}|D\psi|^{p_0}\,d\mu\biggr)^{1/(p_0-1)}\int_{B(x_0,r_0)\setminus \overline{B}(x_0,r)}|DG'|^{p_0}\,d\mu \nonumber \\
& \leq \biggl(\frac{C_4}{C_1}\biggr)^{1/(p_0-1)}(1-\frac{r}{r_0})^{-p_0}\Capc_{p_0}(\overline{B}(x_0,r),B(x_0,r_0))^{1/(p_0-1)} \cdot \nonumber \\
& \int_{B(x_0,r_0)\setminus \overline{B}(x_0,r)}|Dw|^{p_0}\,d\mu \nonumber. 
\end{align}
Let us introduce the function $\xi \in N^{1,p_0}(B(x_0,R))$ defined by
\[
\xi = \left\{\begin{array}{ll}
0, & \textrm{ in } \Om\setminus B(x_0,R), \\
\max\{w,0\}, & \textrm{ in } B(x_0,R)\setminus B(x_0,r_0), \\
w, & \textrm{ in } B(x_0,r_0)\setminus B(x_0,r) \\
\min\{w,M_w(r)\}, & \textrm{ in } B(x_0,r).
\end{array} \right.
\]
Observe that we have $\xi = M_w(r)$ in a neighborhood of $x_0$. 
Let
\[ 
I = \{x\in B(x_0,R):\ \xi (x)=w(x)\}. 
\]
Since $|D\xi|=|Dw|=|DG'|$ on $I$, and $|D\xi| = 0$ on 
$B(x_0,R)\setminus I$, from \eqref{eq:Ku} we have
\begin{align*}
& \int_{B(x_0,r_0)\setminus \overline{B}(x_0,r)}|Dw|^{p_0}\,d\mu \leq \int_I
|Dw|^{p_0-2}Dw\cdot Dw\, d\mu \\
& = \int_I|Dw|^{p_0-2}Dw\cdot D\xi\, d\mu = \int_{B(x_0,R)}|Dw|^{p_0-2}Dw\cdot D\xi\, d\mu \\
& = K(G')M_w(r).
\end{align*}
By plugging this in \eqref{eq:upper}, we finally conclude that
\begin{multline*}
M(r) \geq \biggl(\frac{C_1}{C_4}\biggr)^{1/(p_0-1)}(1-\frac{r}{r_0})^{p_0} \\
\cdot\biggl(\frac{K(G')}{\Capc_{p_0}(\overline{B}(x_0,r),B(x_0,r_0))}\biggr)^{1/(p_0-1)} + M(R).
\end{multline*}
This completes the proof.
\end{proof}

\medskip

\begin{remark}
By obivious modifications, the preceding argument holds in the case 
$p_0=Q(x_0)$ as well.
\end{remark}

\medskip

\begin{remark}
Observe that assuming only conditions 1--3 in Definition~\ref{def:G},
factor $K(G')$ comes up in the above estimates as opposed to the estimates in 
Theorem~\ref{thm:blowupI}.
\end{remark}

\medskip

We have the following result on the local behavior of a Cheeger singular
function near a singularity. The proof of this result is similar to 
that of Theorem~\ref{thm:localbI}, thus, we omit the proof.

\medskip

\begin{thm} \label{thm:localb}
Let $\Om$ be a relatively compact domain in $X$, and $x_0 \in \Om$.
If $G'$ is a Cheeger singular function,
i.e. $G'$ satisfies 
conditions 1--3 in Definition~\ref{def:G}, with singularity at $x_0$, then 
there exist positive constants $C_1,C_2$ and 
$R_0$ such that for any $0<r<\frac{R_0}{2}$ and $x \in B(x_0,r)$ we have
\begin{align*}
C_1\biggl(\frac{d(x,x_0)^{p_0}}{\mu(B(x_0,d(x,x_0)))}\biggr)^{1/(p_0-1)}
& \leq G'(x) \\
& \leq C_2\biggl(\frac{d(x,x_0)^{p_0}}{\mu(B(x_0,d(x,x_0)))}\biggr)^{1/(p_0-1)},
\end{align*}
when $1<p_0<Q(x_0)$, whereas
\[
C_1 \log\left(\frac{R_0}{d(x,x_0)}\right) \leq G'(x) 
\leq C_2\log\left(\frac{R_0}{d(x,x_0)}\right),
\]
when $p_0=Q(x_0)$. 
Here the constants $C_1$ and $C_2$ depend on $K(G')$, $p_0$, $x_0$, 
and the local parameters of $\Om$, and $R_0$ depends only on $\Om$.
\end{thm}

\medskip

The following lemma is well-known and we omit the proof.

\medskip

\begin{lemma} \label{poteq}
Let $K$ be a closed subset of a relatively compact domain $\Om$, and
let $u$ be the \p$_0$-potential of $K$ with respect to $\Om$. Then for all
$0\leq\alpha < \beta \leq 1$ one has
\[
\Capc_{p_0}(\Om^{\beta}, \Om_{\alpha})=\frac{\Capc_{p_0}(K,\Om)}
{(\beta-\alpha)^{p_0-1}}.
\]
\end{lemma}

We close this paper with the following observation. The proof of 
Proposition~\ref{prop:scaling} is similar to the proof of Lemma 3.16 in
Holopainen~\cite{Holopainen}, but we present
it here for completeness.

\medskip

\begin{prop} \label{prop:scaling}
Let $G'$ be a Cheeger singular function, i.e. $G'$ satisfies 
conditions 1--3 in Definition~\ref{def:G}. Then
\[
G = K(G')^{-1/(p_0-1)}G'
\] 
is a (Cheeger) Green's function with equality in condition 4 in 
Definition~\ref{def:G}.
\end{prop}

\medskip

\begin{proof}
Observing that the function $\phii = \min\{G',1\}$ 
can be used in \eqref{eq:Ku}, 
and since $G'$ is the \p$_0$-potential of
the set $\{x\in\Om: G'\geq 1\}$ with respect to $\Om$, 
we obtain
\begin{equation} \label{eq:G'pot}
\Capc_{p_0}(\{x\in\Om:\ G'(x)\geq 1\},\Om) = K(G').
\end{equation}
Let $0\leq \alpha < \beta$ and suppose first that $\beta\leq 
K(G')^{-1/(p_0-1)}$. Then one has
\begin{align*}
& \Capc_{p_0}(\{x\in\Om: G(x)\geq \beta\}, \{x\in\Om:\ G(x)>\alpha\}) \\ 
& = \Capc_{p_0}(\{x\in\Om:\ G'(x) \geq \beta K(G')^{1/(p_0-1)}\}, \\
& \hspace{1.5cm}\{x\in\Om:\ G'(x)>\alpha K(G')^{1/(p_0-1)}\}) \\ 
& = (\beta-\alpha)^{1-p_0}K(G')^{-1}
\Capc_{p_0}(\{x\in\Om: G'(x)\geq 1\}, \Om) \\
& = (\beta-\alpha)^{1-p_0}.
\end{align*}
Let then assume that
$K(G')^{-1/(p_0-1)} < \beta$. Equation \eqref{eq:G'pot} implies
that
\begin{align*}
& \frac{\Capc_{p_0}(\{x\in\Om: G(x)\geq \beta\}, \Om)}{(K(G')^{-1/(p_0-1)}/\beta)^{p_0-1}} \\
& = \Capc_{p_0}(\{x\in\Om: \frac{G(x)}{\beta}\geq \frac{K(G')^{-1/(p_0-1)}}{\beta}\}, \Om) \\
& = K(G'),
\end{align*}
from which it follows that
\[
\Capc_{p_0}(\{x\in\Om: G(x)\geq \beta\}, \Om) = \beta^{1-p_0}.
\]
Then one has
\begin{align*}
& \Capc_{p_0}(\{x\in\Om: G(x)\geq \beta\}, \{x\in\Om:\ G(x)>\alpha\}) \\
& = \Capc_{p_0}(\{x\in\Om:\ G(x)/\beta \geq 1\}, 
\{x\in\Om:\ G(x)/\beta>\alpha/\beta\}) \\ 
& = (1-\alpha/\beta)^{1-p_0}\Capc_{p_0}(\{x\in\Om: G(x)\geq \beta\}, \Om) \\
& = (\beta-\alpha)^{1-p_0}.
\end{align*}
This completes the proof.
\end{proof}

\end{document}